\documentstyle[11pt]{article}

\title{A CIRCUIT-THEORETIC ANOMALY RESOLVED BY NONSTANDARD ANALYSIS}
\author{A. H. Zemanian}
\date{}

\begin{document}
\newcommand{\R} {I \kern -4.5pt R}
\newcommand{\N} {I \kern -4.5pt N}
\maketitle
\baselineskip21pt

{\ Abstract --- An anomaly in circuit theory is 
the disappearance of some of the stored energy when two capacitors, 
one charged and the other uncharged, are connected together through 
resistanceless wires.  Nonstandard analysis shows that, when the wires are 
taken to have infinitesimally small but nonzero 
resistance, the energy dissipated 
in the wires equals that substantial amount of energy that 
had disappeared, and that all but an infinitesimal amount of this
dissipation occurs during an infinitesimal 
initial time period.  This provides still another 
but quite simple model of what is in fact a multifaceted physical 
phenomenon. It also exemplifies the efficacy of at least 
one application of nonstandard analysis to circuit theory.\\

Key Words: Capacitive circuits, capacitive energy anomaly, 
nonstandard circuit analysis
} 

\section{Where Did the Energy Go?}

A much-discussed
anomaly concerning the flow of energy in purely capacitive networks
with switches and resistanceless wires is illustrated perhaps in
its simplest form by the circuit of Figure 1.  Here, two 
capacitors of equal capacitance $c>0$ are initially disconnected 
from each other during time $t<0$ because of an open switch.
Let the charge on the left-hand capacitor be of constant value
$q_{0}=q_{1}(t)$ and the charge on the 
right-hand capacitor be $0=q_{2}(t)$ for $t<0$.
Thus, during that time, $v_{0}=q_{0}/c$ is the voltage on the left-hand
capacitor, and the total stored energy in the circuit is $q_{0}v_{0}/2$. 
With the switch being thrown closed at $t=0$, the charge on each capacitor 
becomes $q_{0}/2$, and the total stored energy becomes
$q_{0}v_{0}/8+q_{0}v_{0}/8=q_{0}v_{0}/4$. Energy in the amount 
of $q_{0}v_{0}/4$ has disappeared.  Where did it go?

\section{No Problem When the Circuit Has Resistance}

Assume that the wires have some resistance $r>0$, as shown 
in lumped fashion in Figure 2.  Then, the disappearing 
energy is accounted for
by dissipation in $r$.  Indeed, for $t>0$, the current 
$i_{r}(t)$ in the circuit is 
\begin{equation}
i_{r}(t)\;=\;\frac{q_{0}}{rc}\,e^{-2t/rc},  \label{2.1}
\end{equation}
and the corresponding power $p_{r}(t)\;=\;i(t)^{2}r$ is
\begin{equation}
p_{r}(t)\;=\;\frac{q_{0}v_{0}}{rc}\,e^{-4t/rc}.  \label{2.2}
\end{equation}
Hence, the total energy $E_{r}(0,\infty)$ dissipated in $r$ is
\begin{equation}
E_{r}(0,\infty)\;=\;\int_{0}^{\infty}\,p_{r}(t)\,dt\;=\;
\frac{q_{0}v_{0}}{4}  \label{2.4}
\end{equation}
This is exactly the difference between the initial capacitively
stored energy $q_{0}v_{0}/2$ and the final capacitively stored energy 
$q_{0}v_{0}/4$ occurring in the limit as $t\rightarrow\infty$.

\section{Standard Distribution Theory Does Not Help}

Perhaps we can account for the discrepancy noted in  Section 1
by using standard distribution theory.  After all, as
$r\rightarrow 0$, $i_{r}(t)$ approaches $\delta(t)q_{0}/2$, 
a delta function of size $q_{0}/2$.
Indeed, for each $r>0$, $i_{r}(t)=0$ for $t<0$ and $\int_{0}^{\infty}
i_{r}(t)\,dt\,=\, q_{0}/2$, whereas $\int_{t>T} i_{r}(t)\,dt\rightarrow 0$
for every $T>0$ as $r\rightarrow 0$.  So, why not use this delta function 
in the first place, assume $r$ is not zero (perhaps it is 
a nonzero infinitesimal), 
and calculate the energy dissipation in $r$ as follows:
\[ \int_{-\infty}^{\infty}(i_{r}(t))^{2}\,r\,dt\;=\;\frac{q_{0}^{2}}{4}r
\int_{-\infty}^{\infty} (\delta(t))^{2}\,dt\;=\;? \]
Unfortunately $\delta^{2}$ is not well-defined as a distribution.  
Singular distributions cannot be multiplied together,
according to standard distribution theory.\footnote{However, nonstandard
generalized functions may circumvent this trouble in multiplication;
in this regard, see the Final Remark in Section 7.}
It appears that the
difficulty of Section 1 persists if we try this approach.

\section{More Complicated Models Can Account for the Vanishing Energy}

One can set up in the laboratory the circuit of Figure 1 using
real components and can note that the discrepancy between the 
initial and final stored energies truly occurs.  One must conclude
that the ideal circuit of Figure 1 with perfect elements is 
too simplified to account for the real phenomenon.
The vanishing stored energy might be accounted for in several ways.
There is indeed some resistance in any wire, which will dissipate energy.
The dielectric within the capacitors has nonlinear resistance 
as well, producing more dissipation.  The current produces a magnetic field
and thereby inductance, which will affect the transient state and
thereby the resistive power dissipation over time.
There may be arcing at the switch as it is being closed, 
producing thereby more heat.  Moreover, the varying electric fields in
the capacitors produce magnetic fields and thereby radiation.
A search on the internet\footnote{For instance, 
search ``capacitance energy loss'' in google or yahoo.} 
yields a number of references discussing all this, some examples of which are
\cite{j-s}, \cite{a-w}, \cite{fo}, \cite{pa}.  
Some of that discussion occurs as internet ``chat.''

So, a variety of different models are suggested by this 
multifaceted physical phenomenon.  Which model is 
preferred depends upon a compromise between 
accuracy and simplicity. 
The objective of this brief note is to suggest one more, but quite simple,
model.  It uses nonstandard analysis to extend the $rc$ circuit of 
Figure 2 to the case where $r$ is a positive infinitesimal, a quantity that 
is less than any real positive number but is not negative, and thus is 
effectively zero from the perspective of standard analysis.
Another purpose of the note is to demonstrate the efficacy 
of nonstandard analysis in circuit theory and perhaps 
encourage its use in other engineering analyses.

\section{How Nonstandard Analysis is Used in This Case}

Up to now, we have tacitly restricted all our variables and 
parameters to real numbers. We will continue to take the capacitance $c$
and the initial charge $q_{0}$ and voltage $v_{0}=q_{0}/c$ as
fixed real positive numbers.  But, we will allow all other quantities to
be nonstandard, namely hyperreal numbers---with real numbers being a 
special case of hyperreal numbers.\footnote{As is conventional, we will simply 
write ``hyperreal'' for ``hyperreal number'' and ``real'' 
for ``real number.''}  The hyperreals comprise the enlargement 
$^{*}\!\R$ of the real line $\R$.  $^{*}\!\R_{+}$ and  $\R_{+}$
denote the positive parts of $^{*}\!\R$ and $\R$.  Each real
$a\in\R$ is contained in a set of hyperreals that are infinitesimally 
close to $a$.  That set is called the halo---or synonymously 
the monad---for $a$.  Moreover, for any two reals $a$ and $b$ where $a< b$,
the halos for $a$ and $b$ do not overlap.  Furthermore, the 
hyperreal line $^{*}\!\R$ extends into unlimited (synonymously, 
infinitely large) hyperreals.  

Here is how these ideas can be used to resolve the anomaly.
Instead of setting $r$ exactly equal to 0, we can let $r$ be a 
positive infinitesimal in Figure 2.  It can then be 
shown that, at each real positive time $t$, 
the current $i_{r}(t)$ and power $p_{r}(t)$
dissipated in $r$ are infinitesimals.  On the other hand, 
it can be shown that, at all sufficiently small positive 
infinitesimal time $t$, the current $i_{r}(t)$ and 
power $p_{r}(t)$ are positive unlimited hyperreals.
It can also be seen that the total energy $E_{r}(0,\infty)$
dissipated in the positive infinitesimal $r$ during the real time 
interval $0\leq t <\infty$ is equal to the real value 
$q_{0}v_{0}/4$, and this is so no matter 
how small this infinitesmal $r$ is chosen.  
No longer does that dissipated energy 
disappear---as it did under standard analysis with $r=0$.
Moreover, it can be seen that the energy $E_{r}(\tau,\infty)\;=\;
\int_{\tau}^{\infty}\, p_{r}(t)\,dt$ dissipated in the 
infinitesimal $r$ during the real time interval $\tau\leq t <\infty$, where 
$\tau$ is any real positive time, is also infinitesimal.
Consequently, we can assert that the energy dissipated
in $r$ during the positive part of the time halo around $t=0$
is infinitesimally close to $q_{0}v_{0}/4$.  Done.

\section{The Details}

To explicate all of this, we invoke a few results from nonstandard analysis.  
There have been many expositions of that theory during the 45 years 
since its inception \cite{ro}.  The book \cite{go} lists 41 such
sources appearing before 1998, and \cite{go} is itself a 
well-written introduction to the subject.  A concise 
listing and explanation of the ideas used herein
can be found in \cite[Appendix A]{ze}.  We shall now present
derivations of the thoughts of the preceding section and will refer 
to \cite[Appendix A]{ze} for certain definitions and results of 
nonstandard analysis.

As above, $\R$ and $^{*}\!\R$ denote the real and hyperreal lines
respectively \cite[Appendix A.5]{ze}, and $\R_{+}$ and $^{*}\!\R_{+}$
denote their positive parts.  Also, $\N=\{0,1,2,\ldots\}$ 
denotes the set of natural numbers.
$\langle a_{n}\rangle$ will denote a sequence of real numbers 
$a_{0},a_{1},\ldots,a_{n},\ldots$ indexed by the natural numbers.
One way (but not the only way) of introducing the hyperreals is 
to define them as equivalence classes of sequences
of real numbers.  To specify that equivalence relation,
we choose a nonprincipal ultrafilter $\cal F$.  This is a set of subsets of 
$\N$ satisfying certain conditions \cite[Appendix A.4]{ze}.
Two sequences $\langle a_{n}\rangle$ and $\langle b_{n}\rangle$
are taken to be equivalent if $\{n\in\N: a_{n}=b_{n}\}\in{\cal F}$.
A hyperreal will be denoted by $[a_{n}]$, where the $a_{n}$ 
are the elements of one (i.e., any one) of 
the sequences in the equivalence class.  Thus, if $\langle a_{n}\rangle$
and $\langle b_{n}\rangle$ are equivalent sequences, then 
$[a_{n}]$ and $[b_{n}]$ denote the same hyperreal (i.e., $[a_{n}]=[b_{n}]$).
Each hyperreal is either positive (i.e., $\{n\!: a_{n}>0\}\in {\cal F}$),
or negative (i.e., $\{n\!: a_{n}<0\}\in {\cal F}$),
or 0 (i.e., $\{n\!: a_{n}=0\}\in {\cal F}$).  Only one of these 
conditions will hold \cite[Appendix A.6]{ze}.

An infinitesimal $[a_{n}]$ is a special kind of hyperreal defined
as follows.  If for every $x\in\R_{+}$ we have 
$\{n\!: |a_{n}|< x\}\in{\cal F}$, then $[a_{n}]$ is an infinitesimal.
Similarly, if for every $x\in\R_{+}$ we have 
$\{n\!: |a_{n}|> x\}\in{\cal F}$, then $[a_{n}]$ is called 
an unlimited hyperreal (synonymously, 
an infinitely large hyperreal).  If $[a_{n}]$ is neither infinitesimal nor 
unlimited, it is called appreciable.  The product of an unlimited 
hyperreal and an appreciable hyperreal is an unlimited hyperreal.
We emphasize that these definitions do not depend upon the 
choice of the representative sequence $\langle a_{n}\rangle$
for $[a_{n}]$;  this is consequence of the properties of the 
nonprincipal ultrafilter $\cal F$.

Now, let $r=[r_{n}]$ be a positive infinitesimal.  Consider Figure 2
with $r$ replaced by $r_{n}$.  Then, for 
any natural number $n$ for which $r_{n}>0$ and for any 
$t\in\R_{+}$, the current $i_{r_{n}}(t)$ is given by (\ref{2.1}) and 
the power dissipated in $r_{n}$ is given by (\ref{2.2}), but with $r$ replaced 
by $r_{n}$ of course.

We shall now show that the hyperreal power $p_{r}(t)=[p_{r_{n}}(t)]$
is infinitesimal at each $t\in\R_{+}$.  Remember that $r_{n}$ and 
$p_{r_{n}}(t)$ are both real positive numbers.  Given any $t\in\R_{+}$ and 
given any $\epsilon\in\R_{+}$, 
there exists a $\rho\in\R_{+}$ such that $r_{n}<\rho$ implies 
that $p_{r_{n}}(t)<\epsilon$.  Thus, 
\begin{equation}
\{n\!: p_{r_{n}}(t)<\epsilon\}\;\supseteq\;\{n\!: r_{n}<\rho\}  \label{5.1}
\end{equation}
But, $\{n\!: r_{n}<\rho\}\in{\cal F}$ because $[r_{n}]$ is
infinitesimal.  By the properties of an ultrafilter (in this 
case, any filter) \cite[Appendix A.4]{ze}, it follows that 
$\{n\!:p_{r_{n}}(t)<\epsilon\}\in{\cal F}$, and this 
is so for every $\epsilon\in\R_{+}$. 
So truly, $[p_{r_{n}}(t)]$
is infinitesimal.

In a similar way, it can be shown that the hyperreal
current $i_{r}(t)=[i_{r_{n}}(t)]$ is infinitesimal for each $t\in\R_{+}$.

Furthermore, $p_{r}(t)=[p_{r_{n}}(t)]$ can be shown to be positive
unlimited for every sufficiently small infinitesimal time
$t=[t_{n}]$ by examining sets in $\cal F$ as above.  But,
let us now use a somewhat more concise argument.  For $r=[r_{n}]$
being a positive infinitesimal, $[q_{0}v_{0}/r_{n}c]$
is a positive unlimited hyperreal.  On the other hand,
if the infinitesimal time $t=[t_{n}]\in\,^{*}\!\R_{+}$ is chosen 
so small that $t<r$, then $[e^{-4t_{n}/r_{n}c}]$
is appreciable, being larger than $e^{-4/c}\in\R_{+}$ and
less than 1.  (To do this, just choose the $t_{n}$ such that
$\{n\!: t_{n}/r_{n}<1\}\in{\cal F}$.  Now, $p_{r}(t)=[p_{r_{n}}(t)]$
is the product of the unlimited $[q_{0}v_{0}/r_{n}c]$ and the appreciable
$[e^{-4t_{n}/r_{n}c}]$
and therefore must be positive unlimited.

Again, in the same way we can show that the hyperreal current 
$i_{r}(t)=[i_{r_{n}}(t)]$ is also positive unlimited for all 
sufficiently small infinitesimal $t=[t_{n}]$.  

Consider now the energy $E_{r_{n}}(\tau,\infty)$ dissipated 
in $r_{n}\in\R_{+}$ from real time $t=\tau\in\R_{+}$ to $t=\infty$.
\[ E_{r_{n}}(\tau,\infty) \;=\;\int_{\tau}^{\infty} p_{r_{n}}(t)\,dt
\;=\;\frac{q_{0}v_{0}}{4}\,e^{-\tau/r_{n}c} \]
Again an argument similar to that given for 
$p_{r}(t)=[p_{r_{n}}(t)]$ (see the argument regarding (\ref{5.1}))
shows that the hyperreal energy $E_{r}(\tau,\infty)\,=\,
[E_{r_{n}}(\tau,\infty)]$  dissipated in the infinitesimal
resistor $r=[r_{n}]$ from $\tau$ to $\infty$ is also infinitesimal, whatever 
be the choice of $\tau\in\R_{+}$.  More particularly,
for any hyperreal $t=[t_{n}]>\tau\in\R_{+}$ (possibly, $t$ is positive 
unlimited), the hyperreal energy
\[ E_{r}(\tau,t)\;=\;\left[\int_{\tau}^{t} p_{r_{n}}(x)\,dx\right] \]
remains infinitesimal and less than $E_{r}(\tau,\infty)$
no matter how large the hyperreal time 
$t=[t_{n}]$ is chosen.

Since all these results hold for every $\tau\in\R_{+}$, we can 
interpret them as follows:  {\em The total hyperreal energy 
dissipated in the infinitesimal resistance
$r=[r_{n}]$ during the real time
interval from 0 to $\tau$ is infinitesimally close to (but less than) 
the real value $q_{0}v_{0}/4$, and this is so no matter how small we
choose $\tau\in\R_{+}$.  On the other hand, that infinitesimal difference 
in the energy is dissipated in $r=[r_{n}]$ during all time 
larger than $\tau\in\R_{+}$, again no matter small $\tau$ is.}
Thus, we see that by using an infinitesimal resistor $r$, 
we have a way of accounting for all 
of the initial energy.

\section{Conclusions}

This solves---by means of nonstandard analysis---the ``mystery
of the vanishing energy'' in the following way:
From the perspective of standard mathematics, the infinitesimal 
resistance $r$ is equivalent to zero resistance because it is less than any 
real positive resistance but not negative.  
Nonetheless, all but an infinitesimal part of the missing 
appreciable energy $q_{0}v_{0}/4$ is found as dissipation in the 
infinitesimal $r$ due to the infinitely large, hyperreal, power dissipation 
$p_{r}(t)$ occurring during some initial positive part of the time
halo at $t=0$.  The remaining infinitesimal part of that missing energy 
occurs as infinitesimal power dissipation occurring during all 
subsequent hyperreal time.

{\em A Final Remark}:  The ``solution'' being proposed here may not
be the only nonstandard way of accounting for the anomaly.  
Effectively, our solution
it is saying that, since there is no anomaly when the resistor 
$r$ is positive, there should be no anomaly when $r$ is infinitesimally small
but positive
(a reflection of the transfer principle of nonstandard analysis
\cite[Appendix A.22]{ze}). 
But, one cannot jump to the limit as $r\rightarrow 0$ because 
the anomaly then jumps into view.  
Perhaps there are other nonstandard solutions available.
For instance, perhaps the newer theories concerning 
the multiplication of distributions coupled with nonstandard analysis 
and that justify the square $\delta^{2}$ of the delta function $\delta$
multiplied by an infinitesimal $r$
may be used to account for the vanishing energy 
at the instant $t=0$.  See, for example, \cite[Section 23]{ob}.

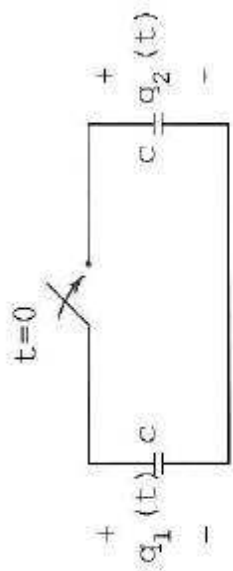
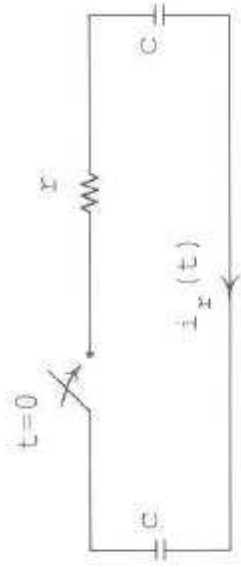
\end{document}